\documentclass[12pt]{article}%[12pt]{amsart}
\usepackage{amssymb,latexsym}
\usepackage{graphicx}
\usepackage{fancyhdr}
\usepackage[colorlinks=true, linkcolor=blue, citecolor=blue,urlcolor=blue, breaklinks=true]{hyperref}
 \usepackage{lineno}
\usepackage{sectsty}
\usepackage{titlesec}
\titleformat{\section}[block]{\Large\bfseries\filcenter}{\thesection}{1em}{}
\begin{document}
\footnotesize
\noindent\framebox[1.02\width]{International Journal of Applied Mathematics 2014; 27 (6), 525-547}
\normalsize
\sectionfont{\large}
\subsectionfont{\normalsize}
\begin{center}
\textbf{MODELING AND NUMERICAL SIMULATIONS OF SINGLE SPECIES DISPERSAL IN SYMMETRICAL DOMAINS} 
%of their free-living stages}Evolution of infectiousness and fitness consequences of
\end{center}
	
%Modeling and analysis of an age-dependent symmetric population model with nonlocality and delay
%\title[Modeling and Numerical Simulation of Population growth] {Modeling and Numerical Simulation of Population growth with delayed nonlocal reaction in a circular domain}\maketitle
\begin{center}
Majid Bani-Yaghoub$^{1}$  Guangming Yao$^{2}$ and Aaron Reed$^{3}$
\end{center}
\begin{table}[h]%\begin{center}
\begin{center}
\begin{tabular}{c}
\noindent $^{1}$ Department of Mathematics and Statistics,\\
                University of Missouri-Kansas City,\\
                Kansas City, Missouri  64110, USA\\
								e-mail: baniyaghoubm@umkc.edu\\
								\\
								\noindent $^{2}$ Department of Mathematics,\\
 Clarkson University,\\
Potsdam, NY, 13699-5815, USA\\
e-mail: gyao@clarkson.edu\\\\
										
\noindent $^{3}$ School of Biological Sciences,\\
University of Missouri-Kansas City,\\
Kansas City, Missouri 64110, USA\\
e-mail: ReedAW@umkc.edu
\end{tabular}%\end{center}
\end{center}
\end{table}
\normalsize
%\vspace{.2cm}
%\author{David E. Amundsen}
%\address{School of Mathematics and Statistics\\
% Carleton University\\
% Ottawa,Ontario, K1S-5B6, Canada}
		\begin{center}
\noindent\textbf{Abstract}
\end{center}
We  develop a class of nonlocal delay  Reaction-Diffusion (RD) models in a circular domain. Previous modeling efforts include RD population models with respect to one-dimensional unbounded domain, unbounded strip and rectangular spatial domain. However, the importance of an RD model in a symmetrical domain lies in the increasing number of empirical studies conducted with respect to symmetrical natural habitats of single species. Assuming that the single species has no directional preference to spread in the symmetrical domain, the RD model is reduced to an equation with no angular dependance. The model can be further reduced by considering the birth function in the form of the Bessel function of the first kind.  We numerically simulate the reduced forms of the nonlocal delay RD model to study the dispersal and growth of behaviors of the single species in a circular domain. Although spatial patterns of population densities are gradually developed, it is numerically shown that the single species population goes extinct in the absence of the birth function or it may converge to a positive equilibrium in the presence of the birth function.\\  
%\noindent \textbf{Math. Subj. Classification:} 37N25 (Dynamical systems in biology), 35R10 (Partial functional-differential equations)\\

\noindent \textbf{Key Words:}  Delay, Reaction-Diffusion, Single Species, Symmetrical Domain\\

%\linenumbers 

\section{Introduction}
\normalsize
\label{intro}
Mathematical modeling of population dynamics has proven to be useful in discovering the relationships between species and their surrounding environment. This includes the study and assessment of spatio-temporal changes in population density and estimation of speed of population dispersal. While various continuous and discrete models have been employed for over a century, recently developed nonlocal delay  Reaction-Diffusion (RD) models have drawn special attention \cite{7 blue, 4 blue, Weng 2008}.  Namely, spatially homogeneous models are equipped with delay, diffusion and integral terms to take into account the  maturation, dispersal and nonlocality of individuals, respectively \cite{1 blue: 9, 6 blue: 8, Thieme 2003}.  The local and global analysis of these models are the current focus of many mathematicians. The present work further develops the age-structured  nonlocal delay RD model of single species proposed by So et al. \cite{7 blue}. To have a better understanding of the nonlocal delay RD model, in the following we briefly explain the modeling procedure initiated by So et. al. \\        
\indent Let $u(t,a,x,y)$ denote the density of the single species at time $t>0$, the age $a\geq 0$ and the spatial position $(x,y)\in\Omega\subseteq\mathbb{R}^{2}$. As described in \cite{6 blue: 48}, dynamics of  the age-structured single species can be formulated  by 
\begin{equation}  \label{Eq:Ch6L1}
\frac{\partial u}{\partial t}+\frac{\partial u}{\partial a} =D(a)\left(\frac{\partial^{2}u}{\partial x^{2}}+\frac{\partial^{2}u}{\partial y^{2}}\right)-d(a)u,
\end{equation}
where $D(a)$ and $d(a)$ are respectively, diffusion and death rates at age $a$. Equation (\ref{Eq:Ch6L1}) describes the spatio-temporal dynamics of single species with respect to age $a$. When there is no age dependence (i.e., when $\partial u/\partial a=0$, $D(a)=D$ and $d(a)=d$ with $D,d>0$), equation (\ref{Eq:Ch6L1}) becomes a linear RD equation that can be directly solved with the method of separation of variables. %Also, in an unbounded domain (i.e. $\Omega=\mathbb{R}^{2}$), with population initially concentrated at origin (i.e. $u(0,x,y)=\delta(x,y)$ where $\delta$ is Dirac delta function).\\
Let $\tau\geq0$ be the maturation time for the single species. Then the total mature population at time $t$ and position $(x,y)$ is given by,
\begin{equation}  \label{Eq:Ch6L2}
w(t,x,y) =\int^{\infty}_{\tau} u(t,a,x,y) da.
\end{equation}
Integrating both sides of (\ref{Eq:Ch6L1}) from $\tau$ to $\infty$, using the assumption $u(t,\infty,x,y)=0$, setting the reproduction density equal to the birth rate (i.e., $ u(t,0,x,y))=b(w(t,x,y))$) and considering the diffusion and death rates to be age independent (i.e., $D(a)=D_{m}$ and $d(a)=d_{m}$ for $a\in[\tau,\infty)$ with $D_{m}$, $d_{m}>0$) we get 
\begin{equation}  \label{Eq:Ch6L3}
\frac{\partial w}{\partial t} =D_{m}\left(\frac{\partial^{2}w}{\partial x^{2}}+\frac{\partial^{2}w}{\partial y^{2}}\right)-d_{m}w+u(t,\tau,x,y).
\end{equation}
Using the procedure outlined in \cite{ Liang, 7 blue, Weng 2008}, $u(t,\tau,x,y)$ can be replaced with an integral term or an infinite series which represents the nonlocality of individuals. When the spatial domain is one-dimensional and unbounded (i.e. $\Omega =\mathbb{R}$), So et. al \cite{7 blue} derived the following population model 
\begin{equation}  \label{Eq:Ch6L4}
\frac{\partial w}{\partial t}=D_{m}\frac{\partial^{2}w}{\partial x^{2}}-d_{m}w+\epsilon\int^{\infty}_{-\infty} b(w(t-\tau, y))f_{\alpha}(x-y)dy,
\end{equation}
where $x\in\mathbb{R}$ and $0<\epsilon\leq 1$. The delay term $\tau>0$ reflects the time required for offspring to become sexually mature. The function $b(w)$ is known as the birth function and reflects reproduction by mature individuals at time $t-\tau$ and any location $y \in \mathbb{R}$. The terms  $D_{m}, D_{I}$ and $d_{m}, d_{I}$ are respectively the dispersal and death rates, where the subscripts $m$ and $I$ respectively correspond to mature and immature population. The kernel function is given by $f_{\alpha}(x)=\frac{1}{\sqrt{4\pi \alpha}}e^{-x^{2}/4\alpha},$ where $\alpha=D_{I}\tau$. Here, $\epsilon$ indicates the total impact of the death rate $d_{I}$ of the immature population which is given by 
\begin{equation}
\label{eqCh3:SoFDE2}
\epsilon=\exp\left\{-\int^{\tau}_{0}d_{I}(a)da\right\}.
\end{equation}
The integral term in (\ref{Eq:Ch6L4}) is a weighted spatial average over the entire spatial domain. Particularly, the integral term is due to the fact that individuals, who are currently at position $x,$ could have been at any location $y \in \mathbb{R}$ at a previous time $t-\tau$. In the present work we will show that the form of the integral term (or equivalently the infinite series for bounded spatial domains) is highly dependent on the the boundary conditions and the shape of the spatial domain.\\
\indent Considering a specific birth function, the traveling wave solution of (\ref{Eq:Ch6L4}) was iteratively constructed in \cite{7 blue}. Later Liang and Wu \cite{3 black} extended the model by adding the advection term $B_{m}\partial w/\partial x$ to the right-hand side of (\ref{Eq:Ch6L4}). Using certain parameter values and birth functions they numerically studied the behavior of traveling wave solutions. Namely, they demonstrated the formation of single and multi-hump wave solutions when the monotonicity condition is violated. %There are also other models similar to equation (\ref{Eq:Ch6L4}) developed by the same methodology. For instance, Al-Omari and Gourley \cite{2 green: 1} generalized the work by Aiello and Freedman \cite{Aiello-Freed} and derived the following model of single species
%\begin{equation}  \label{Eq:Ch6L8}
%\frac{\partial w}{\partial t}=D_{m}\frac{\partial^{2}w}{\partial x^{2}}-d_{m}w+\epsilon\int^{\tau}_{0} G(x,y,s)f(s)e^{-\gamma s} b(w(y,t-s))dy ds,
%\end{equation}
%where $\Omega\subset\mathbb{R}$ is a bounded domain, $f$ is a probability function satisfying $\int^{\tau}_{0}f(s)ds=1$ and $G$ is a kernel derived from solving the heat equation and satisfies $\int_{\Omega}G(x,y,t)dx =\int_{\Omega}G(x,y,t)dy=1$. Moreover, 
Bani-Yaghoub and Amundsen  \cite{bani} showed that a monotonic traveling wavefront of  model (\ref{Eq:Ch6L4})  may become oscillatory when the immature over mature diffusion ratio $D_I/D_m$ is greater than a critical value and the slope of the birth function $b(w)$ at the nontrivial equilibrium is negative (see proposition 2 and section 4 of \cite{bani}). Moreover, Bani-Yaghoub et al.  \cite{bani1} numerically investigated the stability and convergence of solutions associated with model (\ref{Eq:Ch6L4}). They showed that the solution of the initial value problem corresponding to model (\ref{Eq:Ch6L4})  may converge to the corresponding stationary pulse and stationary front.  Ou and Wu \cite{1 blue} developed a general system of RD equations in an m-dimensional domain which embodies a large number of models including (\ref{Eq:Ch6L4}). %Nevertheless their work is limited to the case that delay $\tau\geq0$ is sufficiently small and the existence results are not constructive. 
They showed that for $\tau>0$ sufficiently small, the traveling wavefront exists only if it exists for $\tau=0$. Hence, small maturation time delays $\tau$ are harmless and the traveling wavefronts of the reduced system persist when time lag $\tau$ is increased from zero.  A few other works such as \cite{4 blue} include advection and lift the constraint $\tau\geq0$ being sufficiently small. Nevertheless, they impose other constraints on the kernel function.\\% Nevertheless, due to the fact that the form of the function $f$ in (\ref{Eq:Ch6L8}) is not specified, the outcomes remain at the level of existence and uniqueness and the actual traveling wave solution is not constructed.\\monotonic traveling wave solution of (\ref{Eq:Ch6L4}) becomes oscillatory when the stability of the positive equilibrium is lost.In particular, by considering specific birth functions the convergence of the model solutions to the traveling wave solutions can be numerically investigated \cite{bani1}. Moreover, 
\indent Although model (\ref{Eq:Ch6L4}) is realistic in many aspects, it is constructed with respect to one-dimensional spatial domain, which makes it less appealing. Recent efforts to overcome the issue of single dimension have resulted in models with rectangular spatial domain \cite{Weng 2008} or unbounded strip \cite{Liang}. Nevertheless, dispersal of many single species is radial \cite{Fehmi, Gomes, J} and there is a need to develop models according to symmetrical spatial domains. The present work is an attempt to fill this gap and to capture the spatio-temporal dynamics of single species in circular domains. To have a better understanding of the impact of the spatial domain on the model formulation, in the following we briefly describe two nonlocal delay RD models with different two-dimensional spatial domains.\\
 %Regarding the model %derivation for the case that the spatial domain is two-dimensional, 
%developments, there are two studies that take advantage of the above procedure and derive models with respect to two-dimensional spatial domains. considers the population growth of a single species living in two-dimensional bounded domain. Namely the model
(a)\textit{ Rectangular spatial domain:} The work by Liang et al. \cite{Liang}  considers $\Omega\subset\mathbb{R}^{2}$ as the rectangle $[0,L_{x}]\times[0,L_{y}]$ and zero flux boundary conditions for $u(t,\tau,x,y)$ in time frames $t\in[s,s+\tau]$. Calculating $u(t,\tau,x,y)$ as a function of $w$, equation (\ref{Eq:Ch6L3}) is changed to
\begin{equation}  \label{Eq:Ch6L5}
\frac{\partial w}{\partial t} =D_{m}\left(\frac{\partial^{2}w}{\partial x^{2}}+\frac{\partial^{2}w}{\partial y^{2}}\right)-d_{m}w+F(w(t-\tau,\cdot),x,y),
\end{equation}
where the function $F$ is given by
\begin{equation}  \label{Eq:Ch6L5a}
\begin{array}{ccl}
\vspace{.2cm}
F(x,y,w(t-r,\cdot)) &=& \displaystyle\frac{\epsilon}{L_{x}L_{y}} \int^{L_{x}}_{0}\int^{L_{y}}_{0}b(w(t-r,z_{x},z_{y})).\\
\vspace{.2cm}
&& \Bigg ( 1+\displaystyle\sum^{\infty}_{n=1}\left[\cos\frac{n\pi(x-z_{x})}{L_{x}} +\cos\frac{n\pi(x+z_{x})}{L_{x}}\right]e^{-\alpha(\frac{n\pi}{L_{x}})^2}\\
\vspace{.2cm}
&&+\displaystyle\sum^{\infty}_{m=1}\left[\cos\frac{m\pi(y-z_{y})}{L_{y}} +\cos\frac{m\pi(y+z_{y})}{L_{y}}\right]e^{-\alpha(\frac{m\pi}{L_{y}})^2}\\
\vspace{.2cm}
&&+ \displaystyle\sum^{\infty}_{n=1} \displaystyle\sum^{\infty}_{m=1}\left[\cos\frac{n\pi(x-z_{x})}{L_{x}} +\cos\frac{n\pi(x+z_{x})}{L_{x}}\right] \times\\
\vspace{.2cm}
&&\displaystyle\left[\cos\frac{m\pi(y-z_{y})}{L_{y}} +\cos\frac{m\pi(y+z_{y})}{L_{y}}\right]e^{-\alpha[(\frac{n\pi}{L_{x}})^2+(\frac{m\pi}{L_{y}})^2]}\left. \Bigg )dz_{x}dz_{y}.
\right.\end{array}
\end{equation}
Furthermore, the work by Liang et al. \cite{Liang} derives similar RD models with respect to zero Dirichlet and zero mixed boundary conditions.  %are considered to obtain $F$. 
With certain birth functions, they studied the numerical solutions of the model,  where asymptotically stable steady states and periodic wave solutions were numerically observed. \\
(b) \textit{Unbounded strip:} The work by Weng et al. \cite{Weng 2008} considers a spatial domain $\Omega\subset\mathbb{R}^{2}$ that is a strip in the form of $\Omega=(-\infty,\infty)\times[0,L]$ with $L>0$. Then following the same procedure as outlined in \cite{Liang} they obtain equation (\ref{Eq:Ch6L5}) with the function $F$ defined by 
\begin{equation}  \label{Eq:Ch6L6}
F(w(t-\tau,\cdot),x,y) =\int_{\mathbb{R}}\int^{L}_{0}\Gamma(\alpha,x,z_{x},y,z_{y}) b(w(t-\tau,z_{x},z_{y}))dz_{x} dz_{y},
\end{equation}
where $\Gamma(t,x,z_{x},y,z_{y}) =\Gamma_{1}(t,x,z_{x})\Gamma_{2}(t,y,z_{y})$, $\Gamma_{2} =\frac{1}{\sqrt{4\pi t}}e^{-\frac{(y-z_{y})^{2}}{4t}}$ and $\Gamma_{1}(t,x,z_{x})$ is the Green's function of the boundary value problem,
\begin{equation}  \label{Eq:Ch6L7}
\left\{  \begin{array}{cccl}
 \frac{\partial W}{\partial t}&=& \frac{\partial^2 W}{\partial x^2} &t>0,x\in(0,L)\\
W_B(t,x)&=&0 &t\geq0,x=0,L.
\end{array} \right.
\end{equation}
The term $W_B(t,x)$ denotes zero flux or zero mixed boundary conditions. Using the theory of asymptotic speed of spread and monotone traveling waves, the nonexistence of traveling waves with wave speed $0<c<c^{*}$ and the existence with $c\geq c^{*}$ are established in \cite{Weng 2008}, where $c^{*}$ is known as the minimal speed.\\
\indent In the present paper we will develop a class of  nonlocal delay RD models in a two-dimensional bounded symmetrical domain. The importance of an RD model with symmetrical domain lies in the increasing number of empirical studies conducted with respect to symmetrical natural habitats of single species. Particularly, these studies are conducted by placing the immature population at the center of two-dimensional disks and observing the spread of population over time. For instance, Gomes and Zuben \cite{Gomes} employed a circular arena for radial dispersion of larvae of the blowfly \textit{Chrysomya albiceps}. It is known that after exhaustion of food sources, larvae begin spreading in search of additional food sources.  Then the natural environment can be simulated under experimental conditions by employing circular arenas with sufficiently large diameters (e.g. 50 cm). Also, Roux et al. \cite{Roux} investigated the behavior of the larval dispersal of \textit{Calliphoridae} flies prior to pupation. The study includes statistical results of the shape of the larval dispersal in southwest France in outdoor experimental conditions. The authors found that the shape of the dispersal is circular and has a concentric distribution around the feeding zone. Moreover, the study finds that the larvae had no preference for dispersal in any direction. Although these studies are conducted for dispersal of the larval population, the use of a circular domain and circular dispersal of larvae indicate the need for developing age-structured  nonlocal delay RD models with respect to circular domains.\\ 
\indent The rest of this paper is organized as follows. In section \ref{sec:RDsym} we develop a class of  age-structured  nonlocal delay RD models in circular domains. In section \ref{sec:RadSymm} possible model reductions and the impact of initial heterogeneity are investigated. In section \ref{num} the numerical simulations of the reduced models are presented.   Finally, in Section \ref{disc} a discussion of the main outcomes of this study is provided.  %modeling an , we demonstrate how under some conditions, the two-dimensional model may be reduced to simpler forms.

\section{Model Development}
\label{sec:RDsym}
Focusing on the population of blowflies, the morphological aspects of the larval \textit{Chrysomya albiceps} have been investigated in a number of studies (see \cite{Carvalho} for a review). In particular, there are three stages (i.e. instars) during the larval development of \textit{Chrysomya albiceps} flies. The cephalopharyngeal skeleton of larva develops during the instars and the full development of the skeleton takes place in the third instar. In our study we consider two age classes, where the first two instars are considered as the first age class and the third instar represents the second age class. In addition, the larva displacement takes place always in the landscape and individuals cannot fly. Thus, it is reasonable to consider the two rather than three-dimensional spatial domain. On the other hand we will show that the choice of circular domain can bring valuable insights into the study of symmetric spatial dispersal of individuals.\\
%Otherwise, it remains at the level of a mathematical exercise that is not well-received by the researchers of the broad and diverse areas of expertise.\\
\indent Following the same procedure outlined in section \ref{intro}, assume that equation (\ref{Eq:Ch6L3}) captures the dynamics of blowflies and the spatial domain $\Omega\subset\mathbb{R}^{2}$ is  a two-dimensional disk centered at the origin with radius $R>0$. Since the domain $\Omega$ is a disk, it is suitable to rewrite  equation (\ref{Eq:Ch6L3}) in polar coordinates. 
\begin{equation}  \label{Eq:Ch6L3b}
\begin{array}{ccl}
\vspace{.2cm}
\displaystyle\frac{\partial w}{\partial t} &=& \displaystyle D_{m}\left(\frac{\partial^{2}w}{\partial r^{2}}+\frac{1}{r}\frac{\partial w}{\partial r}+\frac{1}{r^{2}}\frac{\partial^{2}w}{\partial\theta^{2}}\right)-d_{m}w +u(t,\tau,r,\theta).
\end{array}
\end{equation}
Similar to \cite{ Liang, 7 blue, Weng 2008}, we need to replace $u(t,\tau,r,\theta)$ with an explicit function of  $w(t,x,y)$. For $s\geq0$ fixed, define the functional
$$V^{s}(t,r,\theta)=u(t,t-s,r,\theta)\mbox{ with }s<t\leq s+\tau.$$ 
Considering (\ref{Eq:Ch6L1}) in polar coordinates, it follows that for $s\leq t\leq s+\tau$,
\begin{equation}  \label{Eq:Ch6L9}
\begin{array}{ccl}
\displaystyle \frac{\partial V^{s}}{\partial t}(t,r,\theta) &=& \left.\displaystyle \frac{\partial u}{\partial t}(t,a,r,\theta)\right|_{a=t-s} + \displaystyle\frac{\partial u}{\partial a}(t,a,r,\theta) \Bigg|_{a=t-s}, \\
&&\\
&=& \displaystyle D(t-s)\left(\frac{\partial^{2}V^{s}}{\partial r^{2}}+\frac{1}{r}\frac{\partial V^{s}}{\partial r}+\frac{1}{r^{2}}\frac{\partial^{2}V^{s}}{\partial\theta^{2}}\right)-d(t-s)V^{s}.
\end{array} 
\end{equation} 
But note that (\ref{Eq:Ch6L9}) is a linear RD equation that can be solved using the method of separation of variables. Moreover, in the case that the domain is unbounded, the standard theory of Fourier transforms can be used to obtain the general solution of (\ref{Eq:Ch6L9}) (see \cite{hillen} for example). Since $u(t,0,r,\theta)=b(w(t,r,\theta))$, we have
\begin{equation}  \label{Eq:Ch6L9b}
V^{s}(s,r,\theta) =b(w(s,r,\theta)).
\end{equation}
Before solving the initial boundary value problem (IBVP) related to equation (\ref{Eq:Ch6L9}), it would be beneficial to discuss the possible boundary conditions and their biological meanings. Particularly, the choice of the boundary conditions sets certain biological assumptions and it has a major impact on the model development as follows. The zero Dirichlet boundary condition represents the case in which the region outside the domain is uninhabitable. In other words, individuals die once they diffuse out of the domain (see for example \cite{Murray I:292, Murray II}). This makes sense when for instance, individuals are certain aquatic species in a lake or a pond. Nevertheless, zero Dirichlet boundary condition is not suitable for studying species such as amphibians. The book by Kot  \cite{Murray I:292} considers such a boundary condition as an extremely crude way of capturing spatial heterogeneity. Instead, Gurney and Nisbet \cite{GurNisbet} consider that the spatial domain is unbounded and intrinsic rate of growth decreases with the square of the distance from the center of the range. Their approach results in a type of Schr{\"o}dinger equation (see pages 289-291 of \cite{Murray I:292}). The zero-flux boundary condition is another approach that takes away the in-and-out privileges of the individuals. Namely the individuals never cross the boundaries, although they can live and freely move on the boundaries. This has been used in several studies (see for example \cite{bani2, Murray I:292, Murray II}). Combining the zero-flux and Dirichlet boundary conditions gives rise to mixed boundary conditions, where the flux at each boundary is proportional to the population density. Specifically, the individuals may cross the boundary as long as rate exchange with the outer domain at each location remains proportional to the population density at that location. If a boundary is highly populated then we may expect high-population exchange between the inner and outer domains. \\% In the following we treat the problem with respect to all three cases of  boundary conditions mentioned above.\\
\indent Considering zero Dirichlet boundary condition and initial condition described in (\ref{Eq:Ch6L9b}), we have

\begin{equation}  \label{Eq:Ch6L10}%{Eq:Ch6L11}{Eq:Ch6L12}
\left\{
\begin{array}{ccl}
\displaystyle  \frac{\partial^{2}V^{s}}{\partial t} &=& \displaystyle  D(t-s)\left(\frac{\partial^{2}V^{s}}{\partial r^{2}}+\frac{1}{r}\frac{\partial V^{s}}{\partial r}+\frac{1}{r^{2}}\frac{\partial^{2}V^{s}}{\partial\theta^{2}}\right)-d(t-s)V^{s},\\
\displaystyle  V^{s}(t,R,\theta) &=&0,\\
V^{s}(s,r,\theta) &=&b(w(s,r,\theta)).\\
\end{array} \right.
\end{equation}
The IBVP (\ref{Eq:Ch6L10}) can be solved by the method of separation of variables. Specifically, let $V^{s}(t,r,\theta)=h(s,r,\theta)T(t)$; substituting this into (\ref{Eq:Ch6L10}) and separating terms with $h$ from terms with $T$ we find two ordinary differential equations
\begin{equation}  \label{Eq:Ch6L13}
\frac{T^{'}+d(t-s)T}{D(t-s)T} =\lambda,
\end{equation}
\begin{equation}  \label{Eq:Ch6L14}
h_{rr} +\frac{1}{r}h_{r} +\frac{1}{r^{2}}h_{\theta\theta} =\lambda h,
\end{equation}
where $\lambda$ is the separation constant and $(')$ denotes the derivation of $T$ with respect to $t$; $h_{rr}$, $h_{r}$ and $h_{\theta\theta}$ are the partial derivatives of $h$ with respect to $r$ and $\theta$.\\
By letting $\lambda=-k^{2}$ and solving (\ref{Eq:Ch6L13}) we get to
\begin{equation}  \label{Eq:Ch6L14b}
T(t) =\exp\left(-\int^{t}_{s}(k^{2}D(t-\sigma)+d(t-\sigma))d\sigma\right).
\end{equation}
Letting $h(r,\theta)=\rho(r)\Phi(\theta)$ and separating $\rho$ and $\Phi$ in (\ref{Eq:Ch6L14}), we get that the angular part must satisfy 
\begin{equation}  \label{Eq:Ch6L14c}
\Phi^{''}_{n} =-n^{2}\Phi_{n},
\end{equation}
which has the solution
\begin{equation}  \label{Eq:Ch6L15}
\Phi_{n}(\theta) =A_{n}\cos n\theta +B_{n}\sin n\theta,
\end{equation} where $n$ is an integer. The radial equation is
\begin{equation}  \label{Eq:Ch6L16}
r^{2}\rho^{''}_{n} +r\rho^{'}_{n} +(k^{2}r^{2}-n^{2})\rho_{n} =0,
\end{equation} which is the well-studied parametric Bessel equation with solution
\begin{equation}  \label{Eq:Ch6L17}
\rho_{n}(r) =C_{n}J_{n}(kr) +D_{n}N_{n}(kr),
\end{equation} where $J_{n}(kr)$ and $N_{n}(kr)$ are respectively Bessel and Neumann functions of order $n$ and $C_{n}$ and $D_{n}$ are constants. Nevertheless $N_{n}(kr)$ goes to $-\infty$ as $r\rightarrow0$ and we are only interested in bounded solutions. Hence, we set $D_{n}=0$ and $h(r,\theta)$ is written as a linear combination of $h_{n}(r,\theta),$ where
\begin{equation}  \label{Eq:Ch6L18}
h_{n}(r,\theta) =J_{n}(kr)(A_{n}\cos n\theta +B_{n}\sin n\theta).
\end{equation}
In order to satisfy the boundary condition in (\ref{Eq:Ch6L10}), we must have $h(R,\theta) = 0$. This means that $k$ cannot be an arbitrary constant and must satisfy 
\begin{equation}  \label{Eq:Ch6L19}
J_{n}(kR) =0.
\end{equation}
Let $k_{nj}R$ be the $j$-th zero of $n$-th order Bessel function $J_{n}(x)$. Then in equations (\ref{Eq:Ch6L14b}), (\ref{Eq:Ch6L18}) and (\ref{Eq:Ch6L19})  $k$ must be equal to one of the $k_{nj}$s  and the general solution of (\ref{Eq:Ch6L10}) is a linear combination of all these terms, which is given by
\small
\begin{equation}  \label{Eq:Ch6L20}
V^{s}(t,R,\theta) =\sum^{\infty}_{n=0}\sum^{\infty}_{j=1} J_{n}(k_{nj}r)(a_{nj}\cos n\theta +b_{nj}\sin  n\theta)\exp\left(-\int^{t}_{s}k^{2}_{nj}D(t-\sigma)+d(t-\sigma)d\sigma\right). 
\end{equation}
\normalsize
The coefficients $a_{nj}$ and $b_{nj}$ can be determined with the initial condition in (\ref{Eq:Ch6L10}). Let $D_{I}$ and $d_{I}$ denote respectively, the diffusion and death rates of the immature population.  Define 
\begin{equation}  \label{Eq:Ch6L20b}
\epsilon =\exp\left(-\int^{\tau}_{0}d_{I}(a)da\right),
\end{equation}
\begin{equation}  \label{Eq:Ch6L20c}
\alpha =\int^{\tau}_{0}D_{I}(a)da.
\end{equation}
Note that equation (\ref{Eq:Ch6L14b}) can be rewritten as
\begin{equation}  \label{Eq:Ch6L20d}
T(t) =\exp\left(-\int^{t-s}_{0}(k^{2}D(\gamma)+d(\gamma))d\gamma\right).
\end{equation}
When $s=t-\tau$, substituting (\ref{Eq:Ch6L20b})-(\ref{Eq:Ch6L20d}) into (\ref{Eq:Ch6L20}) we have
\begin{equation}  \label{Eq:Ch6L20e}
V^{t-\tau}(t,R,\theta) =\epsilon\sum^{\infty}_{n=0}\sum^{\infty}_{j=1} J_{n}(k_{nj}r)(a_{nj}\cos n\theta +b_{nj}\sin  n\theta)\exp(-k^{2}_{nj}\alpha). 
\end{equation}
Define
\begin{equation}  \label{Eq:Ch6L21}
F_{n}(r) =\sum^{\infty}_{j=1} a_{nj}J(k_{nj}r),
\end{equation}
and
\begin{equation}  \label{Eq:Ch6L22}
G_{n}(r) =\sum^{\infty}_{j=1} b_{nj}J_{n}(k_{nj}r).
\end{equation}
Then for $s=t-\tau$, using the initial condition in (\ref{Eq:Ch6L10}) we have 
\begin{equation}  \label{Eq:Ch6L23}
\sum^{\infty}_{n=0}F_{n}(r)\cos n\theta +G_{n}(r)\sin n\theta =b(w(t-\tau,r,\theta)).
\end{equation}
Equation (\ref{Eq:Ch6L23}) is in the form of Fourier series and therefore $F_{n}(r)$ and $G_{n}(r)$
are given by,
\begin{equation}  \label{Eq:Ch6L24}
F_{n}(r) =\frac{1}{\pi} \int^{2\pi}_{0} b(w(t-\tau,r,\theta))\cos n\theta d\theta, \mbox{  }n=1,2,\ldots,
\end{equation}
\begin{equation}  \label{Eq:Ch6L25}
F_{0}(r) =\frac{1}{2\pi} \int^{2\pi}_{0} b(w(t-\tau,r,\theta))d\theta, \mbox{  } n=0,
\end{equation}
\begin{equation}  \label{Eq:Ch6L26}
G_{n}(r) =\frac{1}{2\pi} \int^{2\pi}_{0} b(w(t-\tau,r,\theta))\sin n\theta d\theta, \mbox{  } n=1,2,\ldots.
\end{equation}
Substituting (\ref{Eq:Ch6L24}) and (\ref{Eq:Ch6L25}) into (\ref{Eq:Ch6L21}), we have
\begin{equation}  \label{Eq:Ch6L27}
\sum^{\infty}_{j=1} a_{nj}J_{n}(k_{nj}r) =\frac{1}{\pi} \int^{2\pi}_{0} b(w(t-\tau,r,\theta))\cos n\theta d\theta, \mbox{  }  n=1,2,\ldots,
\end{equation}
\begin{equation}  \label{Eq:Ch6L28}
\sum^{\infty}_{j=1} a_{nj}J_{n}(k_{nj}r) =\frac{1}{2\pi} \int^{2\pi}_{0} b(w(t-\tau,r,\theta))d\theta, \mbox{  } n=0.
\end{equation}
Similarly, substituting (\ref{Eq:Ch6L26}) into (\ref{Eq:Ch6L22}),  we get
\begin{equation}  \label{Eq:Ch6L29}
\sum^{\infty}_{j=1} b_{nj}J_{n}(k_{nj}r) =\frac{1}{\pi} \int^{2\pi}_{0} b(w(t-\tau,r,\theta))\sin n\theta d\theta, \mbox{  } n=1,2,\ldots.
\end{equation}
For $n$ fixed, each of the series (\ref{Eq:Ch6L27})-(\ref{Eq:Ch6L29}), is recognized as Fourier-Bessel series.  To find the coefficients $a_{nj}$ and $b_{nj}$, we need to multiply both sides by $rJ_{n}(k_{ni}r)$ and integrate from zero to $R$. Thus, from equation (\ref{Eq:Ch6L27}), we have
\begin{equation}  \label{Eq:Ch6L30}
\int^{R}_{0}rJ_{n}(k_{ni}r)  \sum^{\infty}_{j=1} a_{nj}J_{n}(k_{nj}r) dr =\int^{R}_{0}rJ_{n}(k_{ni}r)  \frac{1}{\pi} \int^{2\pi}_{0} b(w(t-\tau,r,\theta))\cos n\theta d\theta,\\
\end{equation}
  with $n=1,2,\ldots.$ \\But note that the Bessel functions are orthogonal with respect to weight function $r$, i.e.,
\begin{equation}  \label{Eq:Ch6L31}
\int^{R}_{0}rJ_{n}(k_{ni}r)  J_{n}(k_{nj}r) dr =0\mbox{ if }k_{ni}\neq k_{nj}.
\end{equation}
Thus, all terms on the left-hand side of (\ref{Eq:Ch6L30}) are zero except the term with $i=j$. We get that
\begin{equation}  \label{Eq:Ch6L32}
a_{ni}\int^{R}_{0}rJ^{2}_{n}(k_{ni}r)dr =\frac{1}{\pi} \int^{R}_{0}\int^{2\pi}_{0} rJ_{n}(k_{ni}r)b(w(t-\tau,r,\theta))\cos n\theta d\theta dr.
\end{equation}
From properties of the Bessel function we have that
\begin{equation}  \label{Eq:Ch6L33}
\int^{R}_{0}rJ^{2}_{n}(k_{ni}r)dr =\frac{1}{2}r^{2}J^{2}_{n+1}(k_{ni}R).
\end{equation}   Therefore,
\begin{equation}  \label{Eq:Ch6L34}
a_{ni} =\frac{2}{\pi R^{2}J^{2}_{n+1}(k_{ni}R)} \int^{R}_{0}\int^{2\pi}_{0} rJ_{n}(k_{ni}r)b(w(t-\tau,r,\theta))\cos n\theta d\theta dr, n=1,2,\ldots.
\end{equation} 
Similarly, applying the same steps to (\ref{Eq:Ch6L28}) and (\ref{Eq:Ch6L29}), we get that
\begin{equation}  \label{Eq:Ch6L35}
a_{0i} =\frac{2}{2\pi R^{2}J^{2}_{1}(k_{0i}R)} \int^{R}_{0}\int^{2\pi}_{0} rJ_{0}(k_{0i}r)b(w(t-\tau,r,\theta))d\theta dr,  \end{equation} 
\begin{equation}  \label{Eq:Ch6L36}
b_{ni} =\frac{2}{\pi R^{2}J^{2}_{n+1}(k_{ni}R)} \int^{R}_{0}\int^{2\pi}_{0} rJ_{n}(k_{ni}r)b(w(t-\tau,r,\theta))\sin n\theta d\theta dr, \mbox{  } n=0,1,2,\ldots.
\end{equation} 
Hence all required elements of the model are determined. Considering that $u(t,\tau,r,\theta)=V^{t-\tau}(t,r,\theta)$, from (\ref{Eq:Ch6L20e}) and (\ref{Eq:Ch6L3b}) we obtain the following  nonlocal delay RD model with initial history function $w_{0}$ and zero Dirichlet boundary condition as follows
\small
\begin{equation}  \label{Eq:Ch6L37} % \label{Eq:Ch6L38} \label{Eq:Ch6L39}
\left\{
\begin{array}{ccl}
\displaystyle  \frac{\partial w (t,r,\theta)}{\partial t}&=&\displaystyle  D_{m}\left(\frac{\partial^{2}w(t,r,\theta)}{\partial r^{2}}+\frac{1}{r}\frac{\partial w(t,r,\theta)}{\partial r}+\frac{1}{r^{2}}\frac{\partial^{2}w(t,r,\theta)}{\partial\theta^{2}}\right)-d_{m}w(t,r,\theta) +\\
&&\displaystyle \epsilon\sum^{\infty}_{n=0}\sum^{\infty}_{i=1} J_{n}(k_{ni}r)(a_{ni}(w(t-\tau,r,\theta))\cos n\theta +b_{ni}(w(t-\tau,r,\theta))\sin n\theta)e^{-k^{2}_{ni}\alpha},\\
w(t,R,\theta) &=& 0\\
w(t,r,\theta) &=&\displaystyle w_{0}(t,r,\theta) \mbox{    for } (r,\theta)\in\Omega, \mbox{ } t\in[-\tau,0]. \\
\end{array} \right.
\end{equation}
\normalsize
where $a_{ni}(w(t-\tau,r,\theta))$ and $b_{ni}(w(t-\tau,r,\theta)))$ are given in (\ref{Eq:Ch6L34})-(\ref{Eq:Ch6L36}), $\alpha$ is defined in (\ref{Eq:Ch6L20c}) and $k_{ni}R$ is the $i$-th zero of $n$-th order Bessel function $J_{n}(x).$ The parameter $\epsilon$ relates to the surviving portion of individuals from birth until they are fully matured. Namely, $0\leq\epsilon <1$ and the portion $1-\epsilon$ of the immature population did not survive and therefore removed from the double sum series in (\ref{Eq:Ch6L37}). \\
%We have used the notations $a_{ni}(w)$ and $b_{ni}(w)$ to emphasize that they are not constants but rather functions of $w$. 
\indent Following the same procedure, we may consider the problem with zero-flux boundary condition and derive a model similar to (\ref{Eq:Ch6L37}). Specifically, in problem (\ref{Eq:Ch6L10}), the boundary condition must be replaced with
\begin{equation}  \label{Eq:Ch6L40}
\frac{\partial V^{s}}{\partial r}(t,R,\theta) =0.  
\end{equation}   Consequently, equation (\ref{Eq:Ch6L19}) is replaced with 
\begin{equation}  \label{Eq:Ch6L41}
\frac{dJ_{n}(kR)}{dr} =0.  \end{equation}
Then $k_{nj}R$ is the $j$-th zero of $n$th order of derivative of the Bessel function (i.e., $dJ_{n}(x)/dx$) and $k$ must be equal to $k_{nj}$ (\ref{Eq:Ch6L14b}) and (\ref{Eq:Ch6L18}). Again, the set of eigenfunctions $\left\{J_{n}(k_{nj}r)\right\}$ form a complete set and they are orthogonal to each other with respect to the weight function in (\ref{Eq:Ch6L30}). Hence the main difference in the model is that in expression (\ref{Eq:Ch6L34})-(\ref{Eq:Ch6L36}), the $\left\{k_{ni}\right\}$ is the set of eigenvalues corresponding to the zero-flux boundary condition (\ref{Eq:Ch6L40}). In particular, the model with zero-flux boundary condition is given by 
\small
\begin{equation}  \label{Eq:Ch6L37b} % \label{Eq:Ch6L38} \label{Eq:Ch6L39}
\left\{
\begin{array}{ccl}
\displaystyle  \frac{\partial w(t,r,\theta)}{\partial t}&=&\displaystyle  D_{m}\left(\frac{\partial^{2}w(t,r,\theta)}{\partial r^{2}}+\frac{1}{r}\frac{\partial w(t,r,\theta)}{\partial r}+\frac{1}{r^{2}}\frac{\partial^{2}w(t,r,\theta)}{\partial\theta^{2}}\right)-d_{m}w(t,r,\theta) +\\
&&\displaystyle \epsilon\sum^{\infty}_{n=0}\sum^{\infty}_{i=1} J_{n}(k_{ni}r)(a_{ni}(w(t-\tau,r,\theta))\cos n\theta +b_{ni}(w(t-\tau,r,\theta))\sin n\theta)e^{-k^{2}_{ni}\alpha},\\
\displaystyle  \frac{\partial w(t,R,\theta) }{\partial r}&=& 0\\
w(t,r,\theta) &=&\displaystyle w_{0}(t,r,\theta) \mbox{    for } (r,\theta)\in\Omega, \mbox{ } t\in[-\tau,0]. \\
\end{array} \right.
\end{equation}
\normalsize
where $k_{ni}R$ is the $i$-th zero of the derivative of $n$-th order Bessel function $J_{n}(x).$   Similarly, a model with nonlocality and delay can be derived with respect to zero mixed boundary condition
\begin{equation}  \label{Eq:Ch6L41b}
A\frac{\partial V^{s}}{\partial r}(t,R,\theta) +BV^{s}(t,R,\theta) =0,  
\end{equation}  where $A$ and $B$ are constants.\\
\indent  Models (\ref{Eq:Ch6L37}) and (\ref{Eq:Ch6L37b}) take into account the angular dependence of $w(t,r,\theta)$ at any location $r$ and time $t$. This means that the population concentrated at origin may have spatial preference in its displacement for search of food or other necessities. However, this is not the case for certain species. As described before, Roux et al. \cite{Roux}, found that there is no preferred direction in spatial movement of the blowfly larvae. Hence we may consider radial symmetry and certain initial conditions to reduce the model into simpler forms. These are discussed in the following section.\\%While the model (\ref{Eq:Ch6L37}) can be numerically solved and existence of traveling wave solutions can be investigated, it is a highly nontrivial problem to find approximations of the solutions. On the other hand,
\section{Model Reduction}
\label{sec:RadSymm}
In the following we assume that population dispersion takes place with radial symmetry but there is no preference at any direction. It follows that the initial condition in (\ref{Eq:Ch6L10}) is independent of $\theta$. Then the solution $V^{s}$ of (\ref{Eq:Ch6L10})  can  also be independent on $\theta$. The ecological interpretation of the initial condition being independent of $\theta$ is that the reproduction of individuals takes place without any angular preference. Specifically, the IBVP (\ref{Eq:Ch6L10}) is reduced to
\begin{equation}  \label{Eq:Ch6L42} % \label{Eq:Ch6L43}\label{Eq:Ch6L44}
\left\{
\begin{array}{ccl}
\displaystyle  \frac{\partial V^{s}}{\partial t} &=& \displaystyle   D(t-s)\left(\frac{\partial^{2}V^{s}}{\partial r^{2}}+\frac{1}{r}\frac{\partial V^{s}}{\partial r}\right)-d(t-s)V^{s},\\
V^{s}(t, R) &=&0, \\ 
V^{s}(s,r) &=&b(w(s,r)). \\
\end{array} \right.
\end{equation}
The zero Dirichlet boundary condition in (\ref{Eq:Ch6L42}) is equivalent to the assumption that the habitat is inhospitable beyond $r=R$. The substitution
\begin{equation}  \label{Eq:Ch6L45}
V^{s}(t,r) =T(t)h(r), 
\end{equation}     reduces equation (\ref{Eq:Ch6L42}) to equation (\ref{Eq:Ch6L13}) and
\begin{equation}  \label{Eq:Ch6L46}
h_{rr}+\frac{1}{r}h_{r} =\lambda h,
\end{equation}  which is a reduced form of equation (\ref{Eq:Ch6L14}). Let $\lambda=-k^{2}$; then (\ref{Eq:Ch6L46}) is rewritten as
\begin{equation}  \label{Eq:Ch6L47}
r^{2}h_{rr} +rh_{r} +r^{2}k^{2}h =0.
\end{equation}
But this is the parametric Bessel equation (\ref{Eq:Ch6L16}) with $n=0$. The solution is given by
\begin{equation}  \label{Eq:Ch6L48}
h(r) =C_{0}J_{0}(kr) +D_{0}N_{0}(kr),
\end{equation}
where $J_{0}$ and $N_{0}$ are respectively Bessel and Neumann functions of order zero. Moreover, $C_{0}$ and $D_{0}$ are arbitrary constants. As indicated before, the Neumann function blows up as $r\rightarrow0$. Specifically,
\begin{equation}  \label{Eq:Ch6L49}
N_{0}(r) \sim\frac{2}{\pi}\ln(\frac{r}{2})\mbox{ as }r\rightarrow0.
\end{equation}
Hence, we let $D_{0}=0$ to obtain a bounded solution for (\ref{Eq:Ch6L42}). It can be shown that
\begin{equation}  \label{Eq:Ch6L50}
J_{0}(r) =\sum^{\infty}_{q=0} \frac{(-1)^{q}}{(q!)^{2}}\left(\frac{r}{2}\right)^{2q}.
\end{equation}   Thus,
\begin{equation}  \label{Eq:Ch6L51}
h(r) =C_{0} \sum^{\infty}_{q=0}\frac{(-1)^{q}}{(q!)^{2}}\left(\frac{kr}{2}\right)^{2q}.
\end{equation}
In order to satisfy the boundary condition in (\ref{Eq:Ch6L42}), we must have $h(R)=0$, then similar to (\ref{Eq:Ch6L19}) we must have
\begin{equation}  \label{Eq:Ch6L52}
J_{0}(kR) =0.
\end{equation}
Let $k_{j}R$ be the $j$-th zero of the Bessel function of order zero; then $k$ must be equal to one of $k_{j}$s. Using (\ref{Eq:Ch6L14b}),   (\ref{Eq:Ch6L45}) and (\ref{Eq:Ch6L48})  the solution of the IBV (\ref{Eq:Ch6L42}) is given by
\begin{equation}  \label{Eq:Ch6L52b}
V^{s}(t,r) =\sum^{\infty}_{j=1} c_{j}J_{0}(k_{j}r)\exp\left(-\int^{t}_{s}k_j^{2}D(t-\sigma)+d(t-\sigma)d\sigma\right).
\end{equation}
Set $s=t-\tau$; using (\ref{Eq:Ch6L20b})-(\ref{Eq:Ch6L20d}) we get to
\begin{equation}  \label{Eq:Ch6L53}
V^{t-\tau}(t, r) =\epsilon\sum^{\infty}_{j=1} c_{j}J_{0}(k_{j}r)\exp(-k^{2}_{j}\alpha). 
\end{equation}
The constant $\alpha$ is defined in (\ref{Eq:Ch6L20c}) and the coefficients $c_{j}$ are determined by the initial condition in (\ref{Eq:Ch6L42}). Namely,
\begin{equation}  \label{Eq:Ch6L54}
b(w(t-\tau,r)) =V^{t-\tau}(t,r).
\end{equation}
Using the fact that (\ref{Eq:Ch6L53}) represents a Fourier-Bessel series, by orthogonality of the Bessel functions and (\ref{Eq:Ch6L33}) we get that,
\begin{equation}  \label{Eq:Ch6L55}
c_{j} =\frac{2}{R^{2}J^{2}_{1}(k_{j}R)} \int^{R}_{0}rJ_{0}(k_{j}r) b(w(t-\tau,r)) dr, \mbox{ } j=1,2,\ldots.
\end{equation}
Hence, the population model of individuals with no directional preference in their spatial dispersal is given by
\begin{equation}  \label{Eq:Ch6L56}
\frac{\partial w}{\partial t} =D_{m}\left(\frac{\partial^{2}w}{\partial r^{2}}+\frac{1}{r}\frac{\partial w}{\partial r}\right)-d_{m}w +\epsilon\sum^{\infty}_{j=1} c_{j}w(t-\tau,r)J_{0}(k_{j}r)\exp(-k^{2}_{j}\alpha),
\end{equation}   with $0\leq r\leq R$, $t\geq0$ and $c_{j}w(t-\tau,r)$ is given in (\ref{Eq:Ch6L55}). The reduced model (\ref{Eq:Ch6L56}) is subject to the initial condition 
\begin{equation}  \label{Eq:Ch6L57}
w(t,r) =w_{0}(t,r), t\in[-\tau,0],
\end{equation}   
and the zero Dirichlet boundary condition 
\begin{equation}  \label{Eq:Ch6L58}
w(t, R) =0.
\end{equation}    %where , $w_{0}$ is the initial history function and (\ref{Eq:Ch6L58}) is the Dirichlet boundary condition.\\
In comparison with model (\ref{Eq:Ch6L37}), equation (\ref{Eq:Ch6L56}) represents a simpler form of population dynamics. Model (\ref{Eq:Ch6L37b}) can also be reduced in the same manner. But in what follows, we can see that under some conditions these models can be substantially simplified. Namely, the final form of the model can be strongly influenced by the choice of the initial condition.  It is known that if the initial condition happens to be in the shape of a particular mode, then the system will vibrate that mode. Let, for instance,  $b(w(r,\theta,s))$ in (\ref{Eq:Ch6L10}) be in the form of
\begin{equation}  \label{Eq:Ch6L59}
b(w(s,r,\theta)) =f(s)J_{1}(k_{2}r)\cos\theta,
\end{equation} where $R=1, k_2 =3.83$ is the second root of $J_0(kR)$ and $f(s)$ is an arbitrary function of $s$. Considering that $s=t-\tau$ is a fixed value, equation (\ref{Eq:Ch6L20e}) is reduced to
\begin{equation}  \label{Eq:Ch6L60}
V^{t-\tau}(t,R,\theta) =\epsilon f(s) J_{1}(k_{2}r) \cos\theta \exp(-k_{2}\alpha).
\end{equation}
Hence, model (\ref{Eq:Ch6L56}) is reduced to
\begin{equation}  \label{Eq:Ch6L61}
\frac{\partial w}{\partial t} =D_{m}\left(\frac{\partial^{2}w}{\partial r^{2}}+\frac{1}{r}\frac{\partial w}{\partial r}+\frac{1}{r^{2}}\frac{\partial^{2}w}{\partial\theta^{2}}\right)-d_{m}w +\epsilon e^{-k_{2}\alpha} f(t-\tau) J_{1}(k_{2}r) \cos\theta.
\end{equation}

\renewcommand{\topfraction}{0.97}
\renewcommand{\textfraction}{0.01}
\begin{figure}[t]
	%\centering
		\includegraphics[width=0.49\textwidth]{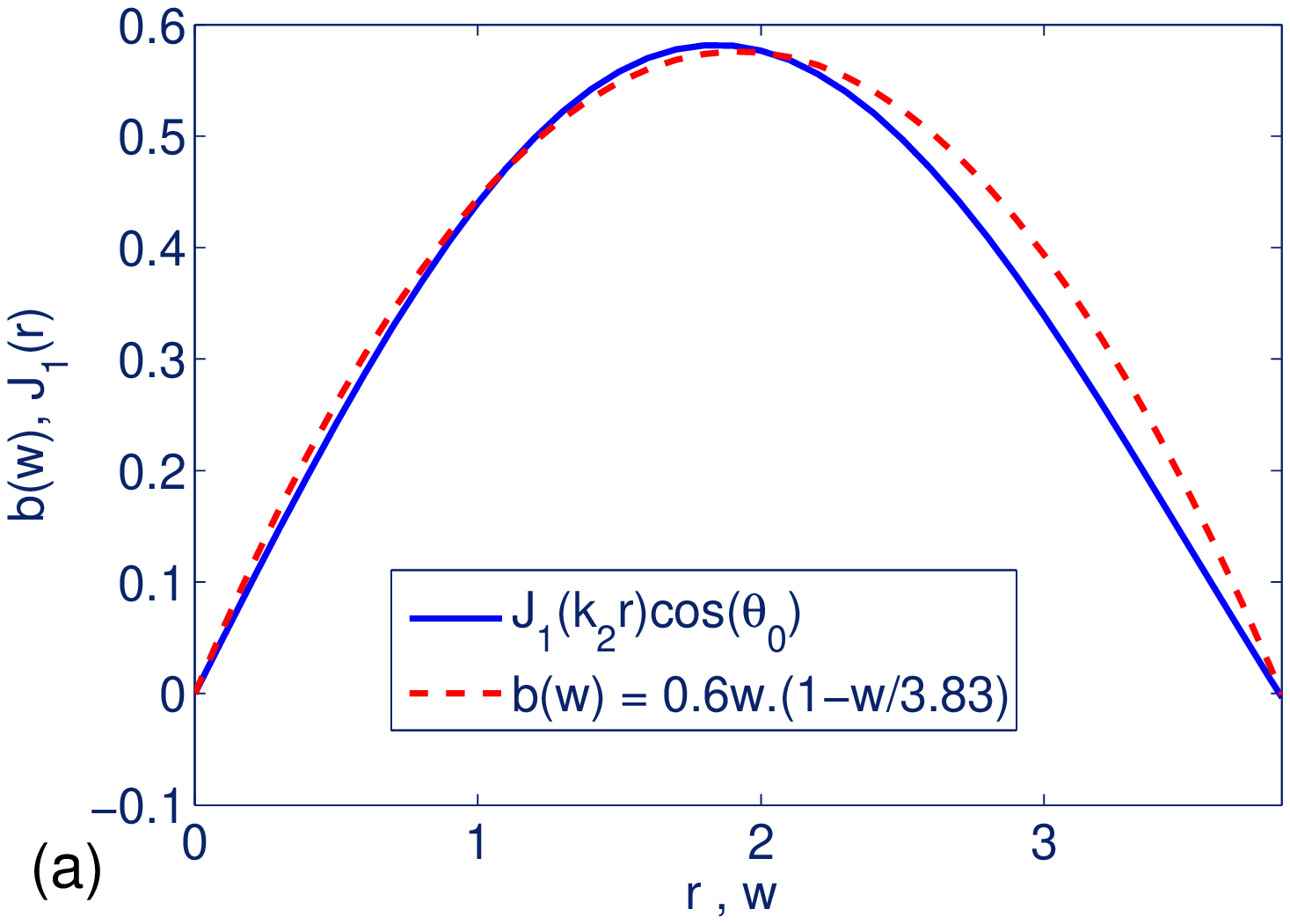}
			\includegraphics[width=0.49\textwidth]{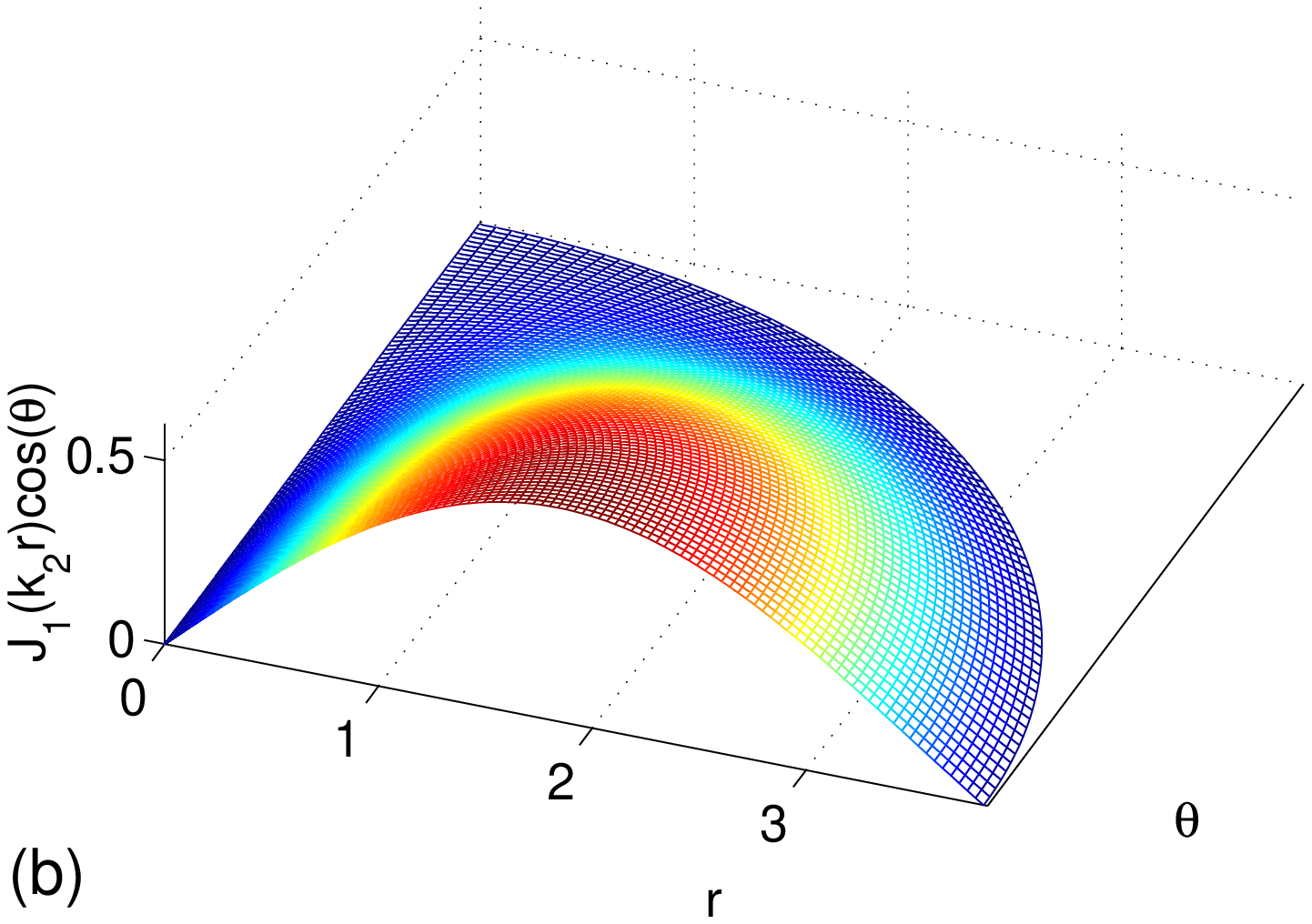}
				\caption{Using equation (\ref{Eq:Ch6L59}) spatially dependent forms of the logistic birth function can be simulated. (a) for $\theta_0 = 0 $ and $f(s)=1$ the initial condition given in equation (\ref{Eq:Ch6L59}) is quite similar to the birth function $b(w) = pw(1-w/k)$. The values of $p$ and $k$ are indicated in the figure. (b) for $0\leq \theta \leq\pi/2  $ the logistic curve preserves its shape while its magnitudes reaches to zero at $\theta =\pi/2$.}
		\label{FigL62} 
\end{figure}

\renewcommand{\topfraction}{0.97}
\renewcommand{\textfraction}{0.01}
\begin{figure}[t]
	\centering
		\includegraphics[width=0.95\textwidth]{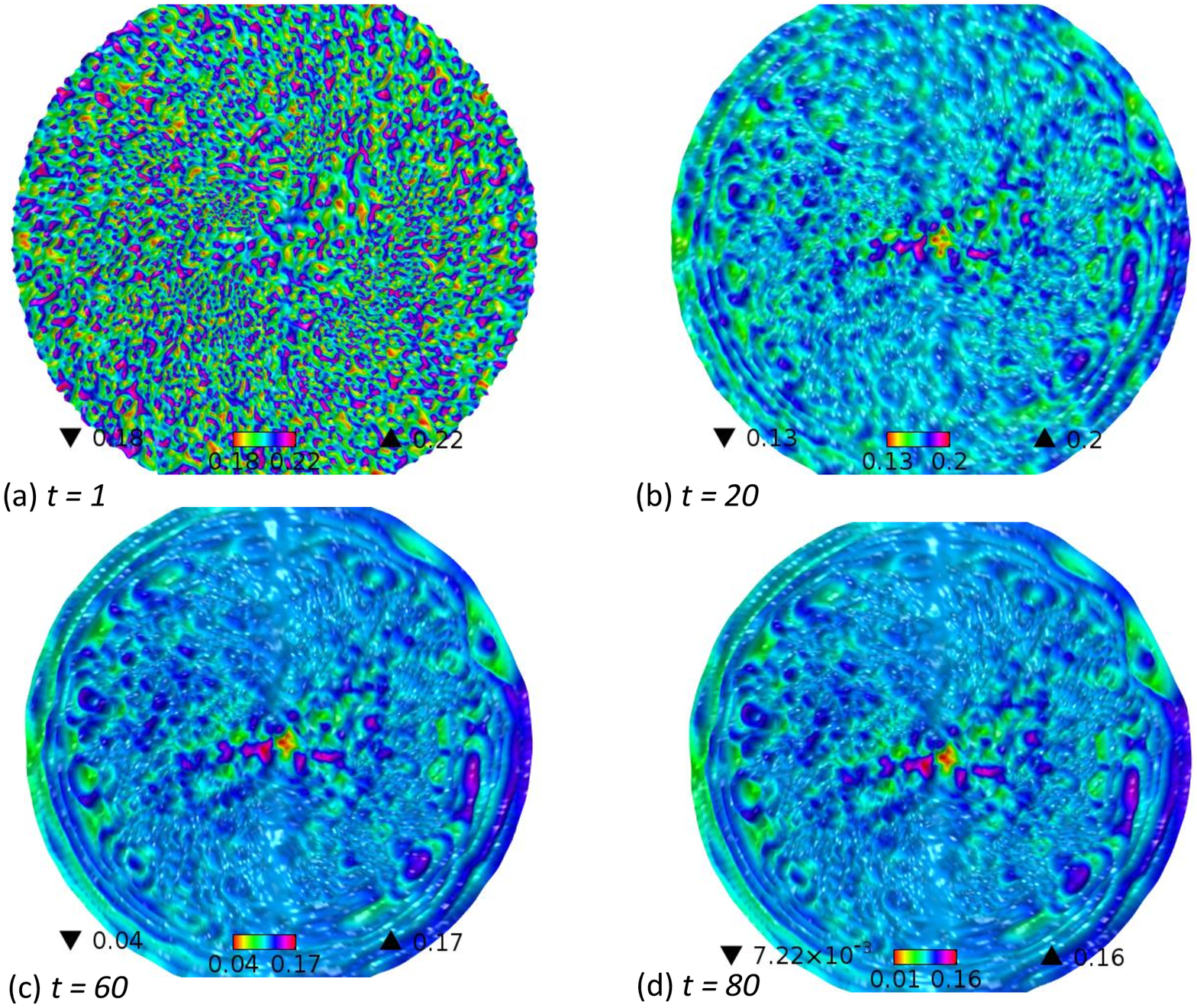}
		\caption{Numerical simulation of the reduced model (\ref{Eq:Ch6L61}). In the absence of  a density dependent birth function $b(w),$ the population goes extinct. Panels (a)-(d) are snapshots of the spatio-temporal variations of the population densities in the xy-plane. Click \href{https://www.youtube.com/watch?v=-g8f5r-06Og&feature=youtu.be}{here} and \href{https://www.youtube.com/watch?v=jUTH8HaJ5n4&feature=youtu.be}{here} to see the transition from (a) to (d) in $xy$-plane and three dimensional cases, respectively.}
		\label{FigL63} 
\end{figure}

The reason for considering the birth function in the form of (\ref{Eq:Ch6L59}) is that for fixed $\theta$ and $s$, the form of the Bessel function $J_{1}(x)$ for $x\in[0, k_{2}R]$ is quite similar to logistic birth function used in several studies \cite{Murray I:292, Murray I}. Figure \ref{FigL62} (a) represents a comparison between the logistic birth function $b(w) = pw(1-w/k)$ and the initial condition given in equation (\ref{Eq:Ch6L59}) for $\theta = 0, f(s)=1, R=1$ and $k=k_2.$ Note that $r$ and $w$ are considered of the same scale. Specifically, the population density increases as we move away from the center of the disk. As shown in Figure \ref{FigL62} (b), by letting $\theta$ change from 0 to $\pi/2$ different birth rates are considered in the spatial domain $[0, R]\times[0, \pi/2].$\\
\indent Although model (\ref{Eq:Ch6L61}) has a much simpler form, it can be argued that the general model (\ref{Eq:Ch6L37}) has been oversimplified  and the reduced model  (\ref{Eq:Ch6L61}) does not fully capture the single species dynamics in a symmetrical domain. Specifically, the reproduction is limited to certain regions of the spatial domain and it is not density dependent.  To overcome these issues we may include the general density dependent birth function $b(w)$ and therefore model (\ref{Eq:Ch6L61}) is rewritten  
\begin{equation}  \label{Eq:Ch6L61b}
\frac{\partial w}{\partial t} =D_{m}\left(\frac{\partial^{2}w}{\partial r^{2}}+\frac{1}{r}\frac{\partial w}{\partial r}+\frac{1}{r^{2}}\frac{\partial^{2}w}{\partial\theta^{2}}\right)-d_{m}w +\epsilon e^{-k_{2}\alpha} f(t-\tau) J_{1}(k_{2}r) \cos\theta + b(w).
\end{equation}
In the next section we will numerically solve the reduced models (\ref{Eq:Ch6L61}) and (\ref{Eq:Ch6L61b}) for different sets of parameter values.
 %The next section deals with numerical solutions of the developed models.\\

\begin{figure}[t]
	\centering
		\includegraphics[width=0.95\textwidth]{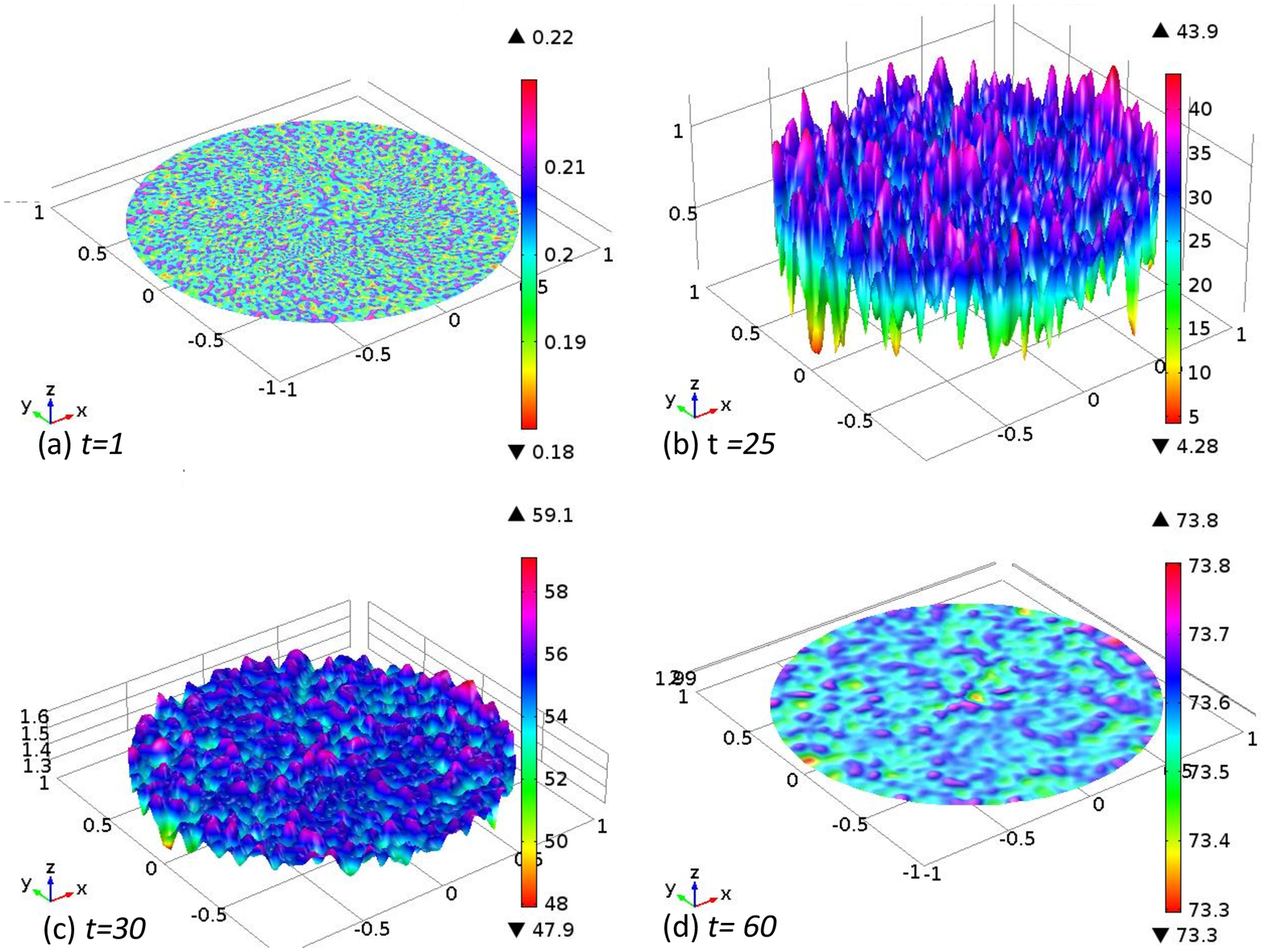}
		\caption{Numerical simulation of the reduced model (\ref{Eq:Ch6L61b}). By considering the birth function $b(w) =.25w^2exp(-0.1w) $ the single species will gradually establish at the positive constant equilibrium. Panels (a)-(d) are snapshots of the spatio-temporal variations of the population densities at specific times.  Click \href{https://www.youtube.com/watch?v=TLUvIoUIe9w&feature=youtu.be}{here} and \href{https://www.youtube.com/watch?v=AeK1HHuCftw}{here} to see the transition from (a) to (d) in $xy$-plane and three dimensional cases, respectively.   }
		\label{FigL64} 
\end{figure}

\section{Numerical Simulations}
\label{num}

We used the COMSOL 4.4 software to solve the reduced models (\ref{Eq:Ch6L61}) and (\ref{Eq:Ch6L61b}). Model (\ref{Eq:Ch6L61}) is a non-homogenous linear RD equation. Due to the presence of  the reaction term $-d_mw,$ it is expected that all solutions $w(r,\theta, t)$ uniformly converge to the trivial solution.  Figure  \ref{FigL63} illustrates the solution of model (\ref{Eq:Ch6L61}) for $t=0$ to $t=80.$  Click \href{https://www.youtube.com/watch?v=-g8f5r-06Og&feature=youtu.be}{here} and \href{https://www.youtube.com/watch?v=jUTH8HaJ5n4&feature=youtu.be}{here} to see the transition from (a) to (d) in $xy$-plane and three dimensional cases, respectively.  The specific parameter values are  $D_m=5, \epsilon = 0.1, d_m=0.01, \alpha = 0.1, R=1, k_2=3.83.$   We considered zero-flux boundary condition and and the initial condition $w_0=0.2+0.02\sin(3x)\cos(2y).$ The color bars in each panel of Figure \ref{FigL63}, show that the population densities are eventually reaching zero. It can be shown that the maximum density is  $5.3 \times 10^{-6}$  when $t > 400$. Hence, as expected, the population goes extinct due to insufficient reproduction and significant mortality of mature population. Model (\ref{Eq:Ch6L61b}) is a non-homogenous nonlinear RD equation, which can be numerically solved. We considered the birth function  $b(w) =.25w^2exp(-0.1w)$ and used the same parameter values as above. Figure  \ref{FigL63} shows model (\ref{Eq:Ch6L61b}) yields the population establishment at the positive constant equilibrium. Click \href{https://www.youtube.com/watch?v=TLUvIoUIe9w&feature=youtu.be}{here} and \href{https://www.youtube.com/watch?v=AeK1HHuCftw}{here} to see the transition from (a) to (d) in $xy$-plane and three dimensional cases, respectively. Using different sets of parameter values, the solutions of model (\ref{Eq:Ch6L61}) converges to the trivial solution. Whereas, the solution of  model (\ref{Eq:Ch6L61b}) has different asymptotic behaviors including the convergence to the positive or trivial equalibria.

\section{Discussion}
\label{disc}
The present work demonstrates that the spatial domain has a great impact in the final form of the derived model. While considering a two-dimensional spatial domain seems to be more realistic, the shape of the spatial domain and the applicable boundary conditions are also important factors that must be carefully dealt with. The work by Weng et al. \cite{Weng 2008} considers an unbounded strip whereas the work by Liang et al. \cite{Liang}  considers a rectangular domain. We believe that model (\ref{Eq:Ch6L5}) with function $F$ specified in (\ref{Eq:Ch6L6}) can also be derived by employing the Smith-Thieme approach \cite{Thieme 2003} for a patchy environment. In particular, the work by So et al. \cite{6 blue: 65} demonstrates how the lattice delay differential equations, representing a population distributed in a line of infinitely many patches, can be extended to the continuous model (\ref{Eq:Ch6L4}) with delay and nonlocality. Replacing the line of infinitely many patches with an  unbounded strip allows us to take into account the spatial movement of individuals within each patch. Then the following the same approach as in (\cite{Gourley3}, pages 5122-5125) the corresponding continuous model is derived. \\
\indent Despite the modeling efforts with respect to unbounded strip and rectangular domains, there is a special need to focus on the symmetrical spacial domains.  The present work is the first step towards developing nonlocal delay RD models with respect to symmetrical domains. It should be noted that the choice of the spatial domain comes from the fact that a number of experimental studies \cite{Fehmi, Gomes, J, Roux} have been conducted in various circular domains. In general we can see that the model derivation is highly dependent on the linear RD equation (\ref{Eq:Ch6L9}), the spatial domain and the boundary conditions. By letting the radius $R$ of the circular domain going to infinity, the spacial domain will be the entire $xy-$ plane. In this case, the same approach outlined in \cite{7 blue} can be used to derive the following nonlocal delay RD model.
\begin{equation}
\label{eqCh3:SoFDE1b}
\frac{\partial w}{\partial t}=D_{m}\left(\frac{\partial^{2}w}{\partial x^{2}}+\frac{\partial^{2}w}{\partial y^{2}}\right)-d_{m}w+\epsilon\int^{\infty}_{-\infty}\int^{\infty}_{-\infty}b(w(z_{x},z_{y},t-\tau))f_{\alpha}(x-z_{x},y-z_{y})dz_{x}dz_{y},
\end{equation}
where $(x,y) \in\mathbb{R}^{2}$, $0<\epsilon\leq 1$, and $w(x,y,t)$ represents the total mature population. The kernel function is given by $f_{\alpha}(x,y)=\frac{1}{\sqrt{4\pi \alpha}}e^{-\frac{x^{2}+y^{2}}{4\alpha}}$ with $\alpha=\tau D_{I}>0$ and $\tau>0$ is the maturation time. \\
\indent  We numerically solved the reduced models (\ref{Eq:Ch6L61}) and (\ref{Eq:Ch6L61b}) for different sets of parameter values and initial conditions. We showed that the solutions related to  model (\ref{Eq:Ch6L61}) converges  to the trivial solution. Whereas the density dependent birth function $b(w)$ considered in model (\ref{Eq:Ch6L61b}) results in convergence of the solutions to the positive equilibrium. Future studies might include numerical simulations of the general models (\ref{Eq:Ch6L37}) and (\ref{Eq:Ch6L37b}). Furthermore, the traveling and stationary wave solutions of these models might bring valuable insights in the studies of single species population dynamics.\\%\\ We believe the correct choice of the birth function $b(w)$ makes the analysis of each model more meaningful in biological and ecological contexts. On the other hand, there have been great similarities in mathematical modeling of single species and those of infectious diseases. There is a potential for employing the same methodology to derive more realistic nonlocal RD models capable of capturing new aspects of the spread of disease in a spatial domain.\\ 

\begin{center}
\noindent\large \textbf{Acknowledgment\\}
\end{center}
\normalsize
  This work was partially supported by University of Missouri-Kansas City start-up fund MOCode \# KCS21.\\

%-----------------------------------------
%-----------------------------------------
%-----------------------------------------
%-----------------------------------------
%-----------------------------------------
%-----------------------------------------
%-----------------------------------------
%-----------------------------------------
%-----------------------------------------
%-----------------------------------------

\end{document}